\documentclass[english]{article}
\usepackage[T1]{fontenc}
\usepackage[latin1]{inputenc}
\usepackage[letterpaper]{geometry}
\usepackage{float}
\usepackage{amsmath}
\usepackage{graphicx}
\usepackage{color}
\usepackage{epstopdf}
\usepackage{adjustbox}
\usepackage{diagbox}
\usepackage{multirow}
\usepackage{array}
\usepackage{mathrsfs}
\usepackage{amsmath, amssymb, amsthm}
\usepackage{subfigure}
\usepackage{tikz}

\graphicspath{{./Figures/}}

\usepackage{amsfonts} 
\usepackage{color}
\definecolor{black}{rgb}{0,0,0}

\definecolor{red}{rgb}{1,0,0}

\definecolor{blue}{rgb}{0,0,1}

\usepackage{multirow}

\makeatother
\usepackage{adjustbox}
\usepackage{diagbox}
\usepackage{babel}
\usepackage{authblk}
\title{A two-grid preconditioner with an adaptive coarse space for flow simulations in highly heterogeneous media }

\author[1]{Yanfang Yang}
\author[2]{Shubin Fu \thanks{Corresponding Author}}
\author[2]{Eric T. Chung}
\affil[1]{School of Mathematics and Information Science, Guangzhou University, Guangzhou, People's Republic of China}
\affil[2]{Department of Mathematics, The Chinese University of Hong Kong, Hong Kong SAR}

\begin{document}
		\maketitle
\begin{abstract}
	In this paper, we consider flow simulation in highly heterogeneous media that has many practical applications in
industry. To enhance mass conservation, we write the elliptic problem in a mixed formulation and introduce a robust two-grid preconditioner to seek the
solution. We first need to transform the indefinite saddle problem to a positive definite problem by preprocessing steps.
 The preconditioner consists of a local smoother and a coarse preconditioner.  For the coarse preconditioner, we design an adaptive
spectral coarse space motivated by the GMsFEM (Generalized Multiscale Finite Element Method). We test our preconditioner for both Darcy flow
and two phase flow and transport simulation in highly heterogeneous porous media. Numerical results show that the proposed preconditioner is highly robust and efficient.
\end{abstract}	

{{\it keywords: } two-grid preconditioner; multiscale;  two-phase flows.}	

\section{Introduction}
	Simulation of fluid flow processes through porous media is important to many subsurface applications, such
	as reservoir simulation, nuclear waste storage and modeling of ground-water contamination.
	For some coupled problems such as the two-phase flow simulation, the most time-consuming part is solving the
	elliptic problem accurately since the media can contain multiple scales and  high-contrast features. Several model
	 reduction techniques such as upscaling \cite{wu2002analysis,SISC_2009} and multiscale techniques \cite{efendiev2006accurate,chung2015mixed,Yang_enriched2018} can alleviate the computational burden.
	 However, the accuracy of the upscaled solution or the multiscale solution can deteriorate with increasing channel correlation length \cite{arbogast2013ms, arbogast2015two}. Moreover, for the coupled flow transport problems,  errors may accumulate with the advancement of time. Therefore, it is necessary to solve the fine-scale flow problem for some cases.  Our
	 goal is to design an effective and  two-grid preconditioner combined with a Krylov accelerator to get the fine-scale solution iteratively. The main feature of the preconditioner is that  it uses the idea of the Generalized Multiscale Finite Element Method (GMsFEM) \cite{efendiev2013generalized,ceh2016adaptive} to form multiscale coarse space. Using the idea of multiscale coarse space to design
	 preconditioner is not new. However, most of the work \cite{galvis2010domain,galvis2010domain2,dolean2012analysis,kim2018bddc,kim2017bddc,klawonn2015feti,mandel2007adaptive,calvo2016adaptive} are devoted to
	 the second order formulation of the elliptic problems.  Here, we focus on the mixed formulation for elliptic problems. The mixed methods are important for many applications, such as flows in porous media, where good approximation to the velocity and mass conservation are required. More accurate approximation of  the velocity can be obtained by using mixed finite element methods since velocity is treated as an independent variable in the method.
	 In \cite{arbogast2015two}, the authors introduce a two-level
	 preconditioner for heterogeneous elliptic problems in mixed formulation
	 with polynomial coarse space.  As we mentioned earlier, we use multiscale coarse space constructed from GMsFEM. 
	
	As we know, the mixed framework results in saddle point problems which make it difficult to get numerical solutions. Therefore, efficient techniques for solving the discretization system are admired for the application of the mixed method.   In the past several decades, many researchers proposed different iterative 
	 methods for the discretized mixed system. In \cite{glowinski1988domain},
	 the authors introduced a nonoverlapping domain decomposition preconditioner. Mathew proposed an overlapping domain decomposition preconditioner in \cite{dd_mixed}. However, the coarse space he
	 used also consists of polynomials which makes the preconditioner not suitable for
	 highly heterogeneous media. Our work here is based on \cite{pre_mixed_hdiv,mg_mixed,dd_mixed}. We adopt the smoothing techniques introduced in \cite{pre_mixed_hdiv} and the preprocessing techniques in \cite{dd_mixed}. We also use more efficient two-grid  method instead of Schwarz method to accelerate the iterative steps.
	 We will incorporate the idea of GMsFEM to design enriched coarse space
	 for the coarse preconditioner. The GMsFEM provides a systemically
	 way to construct coarse space that can capture the major complicated features of the media. The main steps of GMsFEM are first create a rich snapshot
	 space and then select the eigenvectors of carefully  designed
	 local spectral problems corresponding to small eigenvalues. The dimension of the coarse space can be controlled by a pre-defined eigenvalue tolerance.
	 Our preconditioner generally consists of two major components:
	 the local smoother and coarse preconditioner, both parts are very
	 important to the performance of the proposed preconditioner.		
	
	 We test the performance of our preconditioner for both the static
	 Darcy flow simulation and the two-phase flow and transport simulation with a 2-D model and
	  two representative 3-D models. Numerical results show that the proposed preconditioner is highly robust and efficient comparing to other preconditioners that incorporate RT0 and the standard MsFEM space for coarse preconditioner. For the two-phase flow and transport simulation, we
	 only compute the coarse space for the initial permeability field,
	 and keep it fixed with the advancement of the time. This can provide
	 huge computational cost saving, while render good accuracy.
	
	 The rest of the paper is organized as follows. In Section \ref{sec:introd} we first  presents
	 some preliminaries, including grids discretization, the mixed formulation of elliptic problems and its finite element discretization.  A preprocessing step
	 is introduced in Section \ref{sec:pre}, which transforms  saddle point problems to positive-definite problems.  The construction of adaptive coarse space following the GMsFEM is discussed in Section \ref{sec:basis}.
	 Section \ref{sec:2grid} is devoted to describing  the two-grid preconditioner method. In
	 Section \ref{sec:numerical-results}, we presents some representative numerical examples to demonstrate the performance of our preconditioner. A conclusion is drawn in the last section.

\section{Preliminaries}\label{sec:introd}
We consider the following Darcy problem in a mixed formulation:
\begin{equation}\label{eq:orgional_equation}
\begin{split}
\kappa^{-1}v+\nabla p&=0 \quad \text{in } D,\\
\text{div}(v)&=f\quad \text{in } D.
\end{split}
\end{equation}
with the homogeneous Neumann boundary condition $v\cdot n=0$ on $\partial D$, where $\kappa$
is a high-contrast permeability property of the medium, $D$ is the computational domain and $n$ is the unit outward normal
vector of the boundary of $D$, the source function $f$ satisfies $\int_D f=0$. We are mainly interested in computing the
velocity $v$ for the consideration of practical applications such as
reservoir simulation.

To better present our two-grid method, we first introduce the two-scale grid.
Let $\mathcal{T}_H$ be a usual conforming partition of $D$ into quadrilaterals (tetrahedrons for 3D)
$K_i$ with diameter $H_i$ so that $\overline{D}=\cup_{i=1}^N\overline{K}_i$,
where $N$ is the number of coarse blocks.
We call $E_H$ a coarse
face of the coarse element $K_i$ if $E_H = \partial K_i \cap \partial K_j $ or $E_H= \partial K_i \cap \partial{D}$.
Let $\mathcal{E}_H(K_i)$ be the set of all coarse edges (faces) on the boundary of the coarse block $K_i$ and $\mathcal{E}_H=\cup_{i=1}^N\mathcal{E}_H(K_i)$
be the set of all coarse edges (faces). For our adaptive coarse space, the local velocity basis functions are supported on $\omega_i$, which are the two coarse elements that share a common edges (faces), i.e.,
$$\omega_i=\cup \{K\in \mathcal{T}_H: E_i\in \partial K\}, i=1,2,\cdots, N_e,$$ where $N_e$ is the number of coarse faces.
For each coarse block $K_i$, we can define a subdomain $K_i^+$ that
covers $K_i$, therefore $\cup_{i=1}^N K_i^+$ form a non-overlapping decomposition of $D$.
We further partition each each coarse block $K_i$ into a finer mesh with mesh size $h_i$.
Let $\mathcal{T}_h=\cup_{i=1}^N\mathcal{T}_h(K_i)$ be the union of all these partitions, which is a fine mesh partition of the domain $D$.

\begin{figure}[H]
	\centering
	\resizebox{0.7\textwidth}{!}{
		\begin{tikzpicture}[scale=0.7]
		\filldraw[fill=orange, draw=black] (2,6) rectangle (10,14);
		\filldraw[fill=blue, draw=black] (4,4) rectangle (8,12);
		
		\draw[step=4,black, line width=0.8mm] (0,0) grid (12,16);
		\draw[step=1,gray, thin] (0,0) grid (12,16);
		\draw[ultra thick, red](4, 8) -- (8,8);
		\draw [->,dashed, thick](12.4,10) -- (6.5, 8);
		\node at (16.8,10.2)  {\huge $E_{i}$: coarse edge (red)};
		
		\draw [->,dashed, thick](12.4,6) -- (6.5, 7.2);
		\node at (19.2,6.2)  {\huge $\omega_{i}$: Coarse neighborhood (blue)};
		
		\draw [->,dashed, thick](12.4,13) -- (8.5, 12.4);
		\node at (20.8,12.8)  {\huge $K_i^{1,+}$: oversampled coarse
			block (orange)};
		\node at (6.5,10.5)  { \huge $K_i^1$} ;
		\node at (6.5,5.5)  { \huge $K_i^2$} ;
		\end{tikzpicture}
	}
	\caption{ Illustration of a coarse edge $E_i$, and its coarse neighborhood $\omega_{i}$, coarse blocks $K_i^1$ and $K_i^2$ and oversampled coarse
		block $K_i^{1,+}$.}
	\label{grids}
\end{figure}
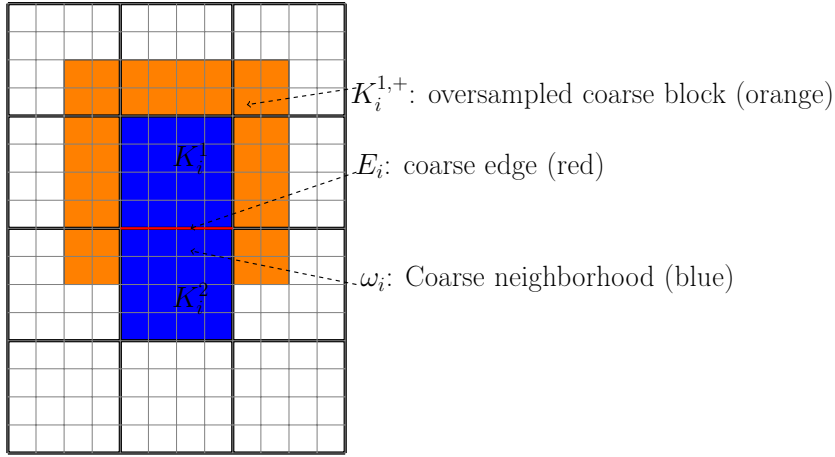

Figure~\ref{grids} gives an example of the constructions of the two-scale grid for the case of 2D. The black lines represent the coarse grid, and the Grey lines represent the fine grid.

We let $V_h\subset H_0(\text{div},D)$ and $Q_h\subset L^2(D)$ to be the lowest-order Raviart-Thomas finite element spaces
with respect to prescribed triangulation $\mathcal{T}_h$
for the approximation of (\ref{eq:orgional_equation}) .
Then,  $v_h\in V_h, p_h\in Q_h$ satisfy
\begin{equation}
\begin{split}
\int_D\kappa^{-1}v_h\cdot w_h-\int_D\text{div}(w_h)p_h&=0, \quad\quad\quad\quad \forall w_h\in V_h^0,\\
\int_D\text{div}(v_h)q_h&=\int_Dfq_h, \quad\quad \forall q_h \in Q_h.
\end{split}
\label{eq:matrixform}
\end{equation}
where $v_h\cdot n=0$ on $\partial D$ and $V_h^0=V_h\cap\{v\in V_h: v\cdot n=0  \text{ on } \partial D\}$.

Let $\phi_i,\cdots,\phi_n$ and $q_1,\cdots,q_m$ be the basis sets for $V_h$ and $Q_h$ respectively,
then assume $v_h=\sum_{i=1}^{n}u_i\phi_i$ and $p_h=\sum_{i=1}^{m}p_iq_i$.
Then above system can be written in terms of matrix representations as
\begin{equation}
\left[\begin{array}{cc}
A_h & B_h^{T} \\
B_h & 0\\
\end{array}\right]\left[\begin{array}{cc}v_{h}\\p_{h}\end{array}\right] =\left[\begin{array}{cc}0\\F_h\end{array}\right]
\label{eq:fine_system}
\end{equation}
where $A_h$ is a symmetric, positive definite matrix with $A_{h,ij}=\int_D\phi_i^T\kappa^{-1}\phi_j$,
$B_h$ is an approximation to the divergence operator with $B_{h,ij}=\int_Dq_i\text{div}\phi_j$,
$F_h$ is a vector with $F_{h,i}=\int_Df_hq_i$

With the coarse and fine grid, we can also define related subspace and operators.
Denote $V_H$ and $Q_H$ be the Raviart-Thomas velocity and pressures space on the coarse grid $\mathcal{T}_H$.
$V_H$ and $Q_H$ are the subspaces of $V_h$ and $Q_h$ respectively.
For each coarse block $K_i$, let $V_i=V_h\cap H_0(\text{div},K_i)$ and $Q_i=Q_h\cap L^2(K_i)$.
Similarly, for each oversampled subdomain $K_i^+$, we denote $V_i^+=V_h\cap H_0(\text{div},K_i^+)$ and $Q_i^+=Q_h\cap L^2(K_i^+)$.
Denote $R_H^T:V_H\times Q_H\to V_h\times Q_h $ as the standard interpolation from coarse
space to the fine space, then $R_H$ is the restriction operator from the fine space to coarse space.
For the subdomains, we let $R_i^T:V_i\times Q_i\to V_h\times Q_h$  be the extension by zero from
$K_i$ to $D$, $R_i$ denotes the restriction map from $V_h\times Q_h$ to $V_i\times Q_i$.
With the same arguments, we can define $R_i^{+T}$ and $R_i^+$.

We can also define some submatrices on local grids and coarse grid in terms of the
restriction and extension operators. We have the coarse grid matrix as
\begin{equation}
L_H\equiv
\left[\begin{array}{cc}
A_H & B_H^{T} \\
B_H & 0\\
\end{array}\right]=R_0\left[\begin{array}{cc}
A_h & B_h^{T} \\
B_h & 0\\
\end{array}\right]R_0^T.
\end{equation}
The coefficient matrix on the subdomain $K_i$ is
\begin{equation}
L_i\equiv
\left[\begin{array}{cc}
A_i & B_i^{T} \\
B_i & 0\\
\end{array}\right]=R_i\left[\begin{array}{cc}
A_h & B_h^{T} \\
B_h & 0\\
\end{array}\right]R_i^T.
\end{equation}
and on the oversampled subdomain $K_i^+$ is
\begin{equation}
L_i^+\equiv
\left[\begin{array}{cc}
A_i^+ & B_i^{+T} \\
B_i^+ & 0\\
\end{array}\right]=R_i^+\left[\begin{array}{cc}
A_h & B_h^{T} \\
B_h & 0\\
\end{array}\right]R_i^{+T}.
\end{equation}

\section{Preprocessing}\label{sec:pre}
We aim to solve the system (\ref{eq:fine_system}) with the preconditioned conjugate gradient (PCG) method, however, the system is not positive
definite, therefore we can not apply PCG directly. We propose an inexpensive preprocessing procedure to transform the problem into a PCG solvable problem.

The idea is to find a discrete flux $\bar{v}_h\in V_h$ such that
\begin{equation*}
\text{div}(\bar{v}_h)=f_h
\label{eq:discrete}
\end{equation*}
where $f_h$ is the $L^2$ projection of $f$ into the space $Q_h$.
Equation (\ref{eq:discrete}) is equivalent to $B_h\bar{v}_h=F_h$ in matrix
form.
To find this discrete flux, we first compute a discrete velocity $\bar{v}_H$
by solving the original equation (\ref{eq:orgional_equation}) on the coarse grid $\mathcal{T}_H$:
\begin{equation}
\left[\begin{array}{cc}\bar{v}_H\\p_H\end{array}\right]=R_H^T
\left[\begin{array}{cc}
A_H & B_H^{T} \\
B_H & 0\\
\end{array}\right]^{-1}R_H \left[\begin{array}{cc}0\\F_h\end{array}\right]
\end{equation}
$\bar{v}_H$ does not necessarily satisfy $\text{div}(\bar{v}_H)=f_h$ on the fine grid $\mathcal{T}_h$. However, $B_h\bar{v}_H-F_h$ has mean value on each coarse block $K_i$ since
$\bar{v}_H$ is obtained by solving the original problem in weak formulation on the coarse grid. Therefore, for  each coarse block $K_i$, we can solve below zero flux boundary
condition subproblem as
\begin{equation}
\left[\begin{array}{cc}{\bar v}_i\\p_i\end{array}\right] =R_i^TL_i^{-1}R_i\left[\begin{array}{cc}-A_h{\bar v}_H\\F_h-B_h{\bar v}_H\end{array}\right]
\label{eq:subproblem}
\end{equation}	
then we let $\bar{v}_h=\bar {v}_H+\bar {v}_1+\cdots+\bar {v}_N$,
it can be verified that  $B_h\bar{v}_h=F_h$.

Then, the solution to (\ref{eq:matrixform}) can be written as
\begin{equation*}
\left[\begin{array}{cc}v_h\\p_h\end{array}\right]=
\left[\begin{array}{cc}\bar{v}_h\\0\end{array}\right]+
\left[\begin{array}{cc}\tilde{v}_h\\p_h\end{array}\right],
\end{equation*}
where $\tilde{v}_h, p_h$ satisfies
\begin{equation}
\left[\begin{array}{cc}
A_h & B_h^{T} \\
B_h & 0\\

\end{array}\right]\left[\begin{array}{cc}{\tilde v}_h\\p_h\end{array}\right] =\left[\begin{array}{cc}-A_h{\bar v}_h\\0\end{array}\right]
\label{eq:reduced_fine_system}
\end{equation}
The above system can be solved with PCG since $\tilde{v}_h$ is divergence free,
and we can solve all subdomain problems in a divergence free space, for more
details we refer \cite{dd_mixed}.

We remark that the preprocessing step is cheap, since it only involves solving
a coarse-grid problem and several subdomain problems and thus can be parallelized without any difficulties.

\section{Adaptive spectral  space}\label{sec:basis}
In this section, we present the construction of adaptive spectral space $V_{H}$ for the velocity in detail. We first need a snapshot space
$V_{\text{snap}}$ from which we can preform model reduction.
The reduction is achieved through a carefully designed local spectral problem.
We select those dominant modes to form  the coarse space.
Notice that we use the terminology of the GMsFEM introduced in \cite{efendiev2013generalized}.
\subsection{Snapshot space}\label{snapshot_space}
In this section, we discuss the formation of the snapshot space  $V_{\text{snap}}$ which consists of basis functions up to the resolution of the fine grid faces on the coarse grid faces.
We construct the local snapshot spaces $V^{i}_{\text{snap}}$ by solving a set of local problems on each coarse neighborhood  $\omega_i$, and then get $V_{\text{snap}}=\oplus_{E_i \in \mathcal{E}_H }V^{i}_{\text{snap}}.$

Let $E_i \in {\cal E}_H$, which can be written as a union of fine-grid faces, i.e., $E_i=\cup_{l=1}^{J_i}e_l$, where $J_i$ is the total number of fine-grid faces on $E_i$ and $e_i$ represents a fine-grid face.
We solve the following problems
\begin{equation}
\begin{split}
\kappa^{-1} v_l^{i}+\nabla p_l^i&=0,\quad\quad \text{in} ~ \omega_i,\\
\text{div}(v_l^{i})&=\alpha_l^i,\quad \quad \text{in}~ \omega_i.
\end{split}
\label{snapshot_eqn}
\end{equation}
subject to the homogeneous Neumann boundary condition $v_l^{i} \cdot n_i=0$ on $\partial \omega_i.$ The above problem is solved separately on each coarse-grid block contained in $\omega_i$, so that the snapshot basis consists of solutions
of local problems with all possible boundary conditions on the face $E_i$ up to the fine-grid resolution.  To solve the equation (\ref{snapshot_eqn}) on $K\subset \omega_i$, an additional boundary condition
$v_l^{i} \cdot n_i=\delta_l^i$ on $E_i$ is used, where $\delta_l^i$ is defined by
\begin{equation}
\delta_l^i=\left\{
\begin{aligned}
1,\quad \text{on} \quad e_l,\\
0, \quad \text{on} \quad E_i \backslash e_l,\\
\end{aligned}
\quad l=1,2, \cdots, J_i,
\right.
\end{equation}
and $n_i$ is a fixed unit-normal vector for $E_i.$ The constant $\alpha_l^i$ in equation (\ref{snapshot_eqn}) is chosen to satisfy the compatible condition $\int_k \alpha_l^i=\int_{\partial K}v_l^i\cdot n_i.$
The solutions of the above local problems form the local snapshot space $V^{i}_{\text{snap}}$, from which we get  $$V_{\text{snap}}=\oplus_{E_i \in \mathcal{E}_H }V^{i}_{\text{snap}}=\text{span}\{v_l^{i}| 1\leq l \leq J_i, 1\leq i \leq N_e\}.$$ Next, we discuss the
derivation of the multiscale space from $V_{\text{snap}}$.

\subsection{Generalized multiscale space}\label{offline_space}
As we mentioned earlier, the snapshot space $V_{\text{snap}}$ is a large space with dimension comparable to the fine grid resolution. To further reduce the space dimension, we will perform a dimension reduction on $V_{\text{snap}}$ to get a smaller space. This reduced space is called
the multiscale space.
The reduction is accomplished by solving a local spectral problem on each coarse grid neighborhood $\omega_i,$ and selecting some dominant modes from the snapshot space $V^{i}_{\text{snap}}$.
The local spectral problem for a coarse face $E_i$ is to find real number $\lambda$ and $v\in V^{i}_{\text{snap}}$ such that
\begin{equation}
a(v,w)=\lambda s(v,w),\quad \forall ~ w \in V^{i}_{\text{snap}},
\label{spectral_1}
\end{equation}
where $a(\cdot,\cdot)$ and $s(\cdot,\cdot)$ are symmetric positive definite bilinear operators on $V^{i}_{\text{snap}}\times V^{i}_{\text{snap}}.$
Specifically, we take
\begin{equation}
\begin{split}
a(v,w)&=\int_{E_i} \kappa^{-1}(v\cdot n_i)(w\cdot n_i),\\
s(v,w)&= \frac 1 {H}\left(\int_{\omega_i} \kappa^{-1}v\cdot w+\int _{\omega_i} \text{div}(v)\text{div}(w)\right),
\end{split}
\label{spectral_eqn}
\end{equation}
for $v, w\in V^{i}_{\text{snap}}$, and $n_i$ is the fixed unit normal vector for $E_i$.

\begin{figure}[H]
	\centering
	\subfigure[local $\kappa$]{
		\includegraphics[width=2.5in]{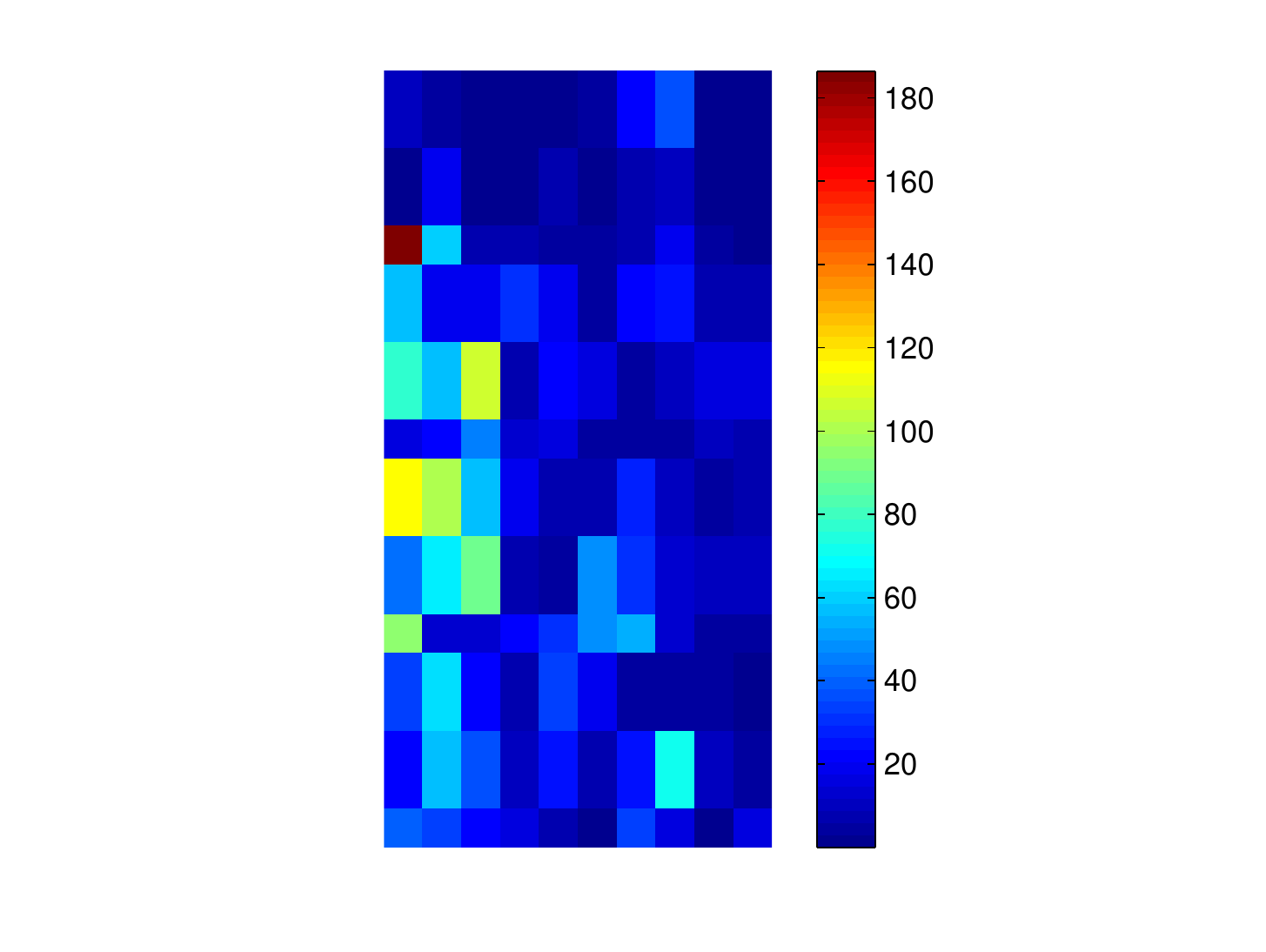}\label{fig:localk} }	
	\subfigure[inverse of the eigenvalues]{
		\includegraphics[width=2.5in]{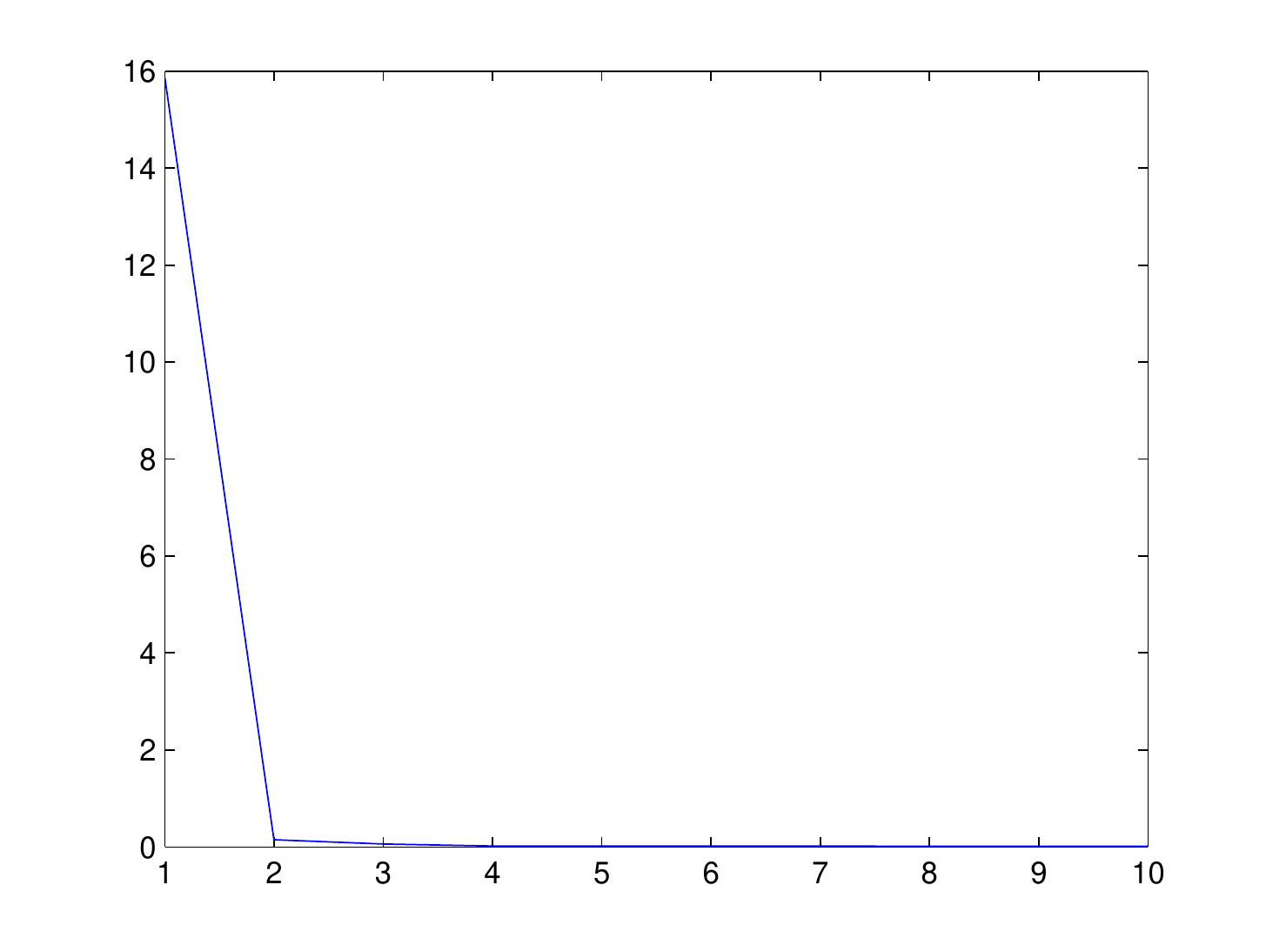}\label{fig:eig} }
		\caption{Local permeability and corresponding decay of the inverse of the eigenvalues}
\end{figure}

\begin{figure}[H]
	\centering
	\subfigure[First spectral basis (x-component)]{
		\includegraphics[width=3in]{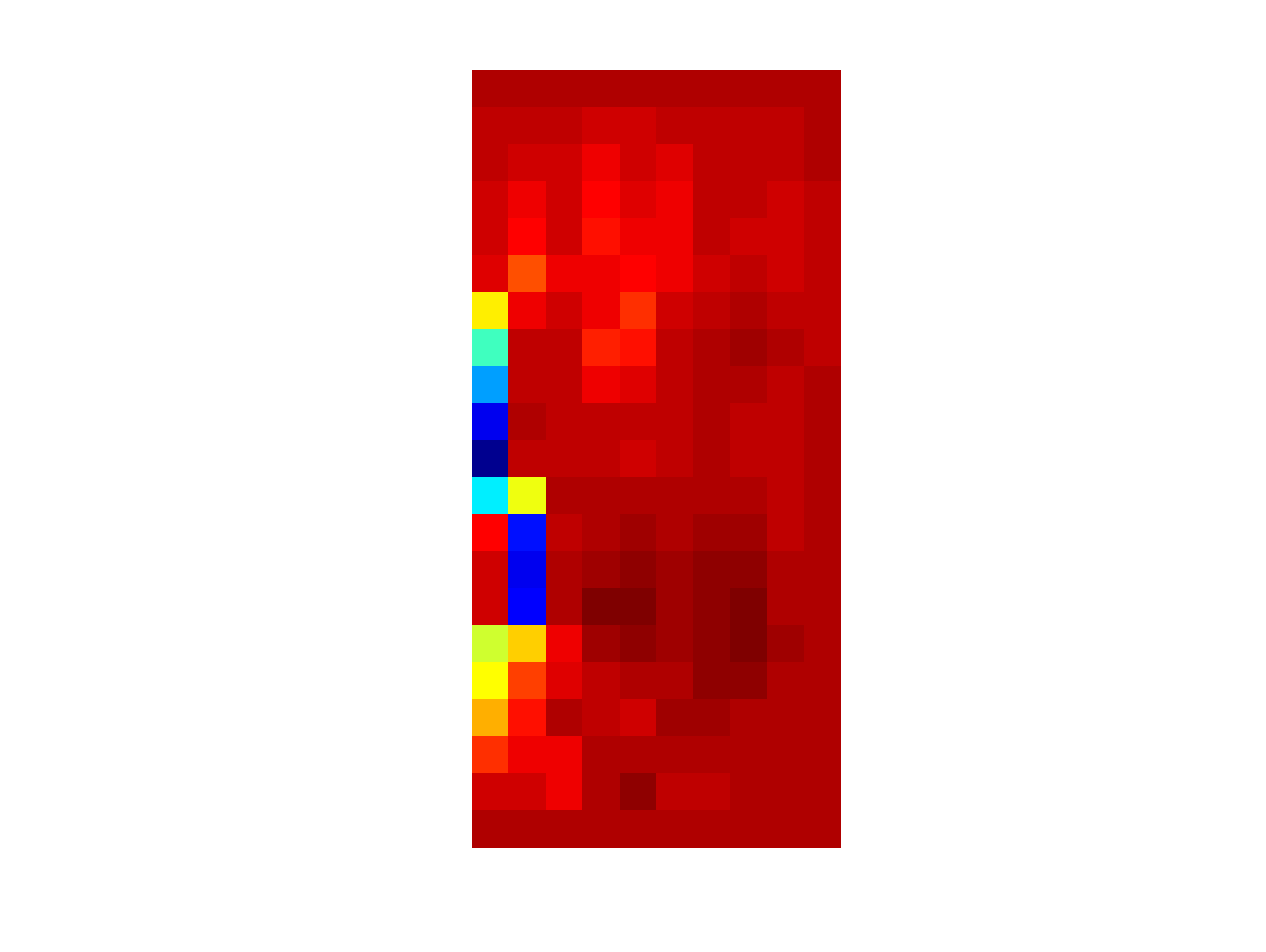}\label{fig:eigb1} }	
	\subfigure[Second spectral basis (x-component)]{
		\includegraphics[width=3in]{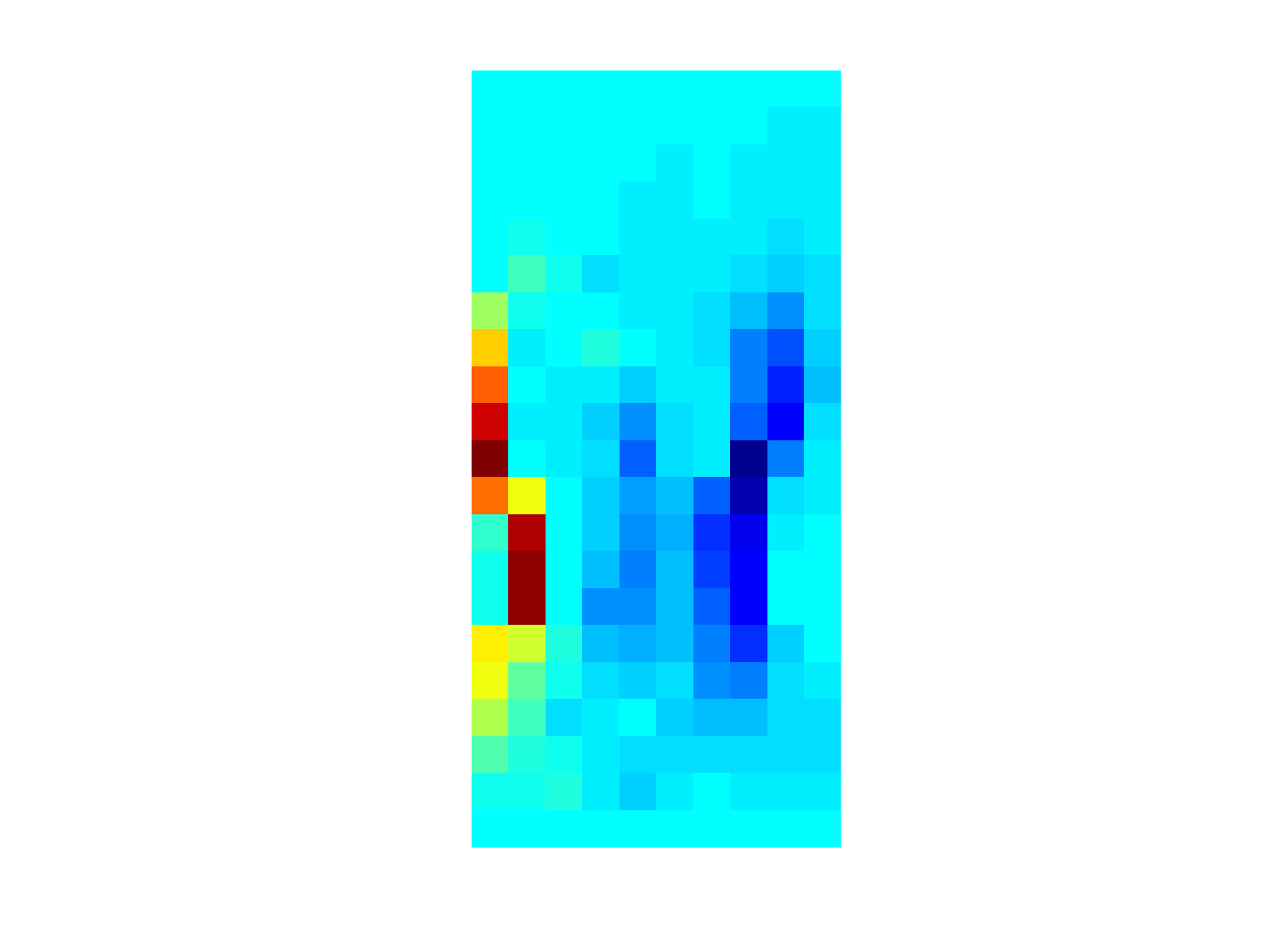}\label{fig:eigb2} }
		\caption{Multiscale basis with GMsFEM (x-component)}
\end{figure}

\begin{figure}[H]
	\centering
\includegraphics[width=3in]{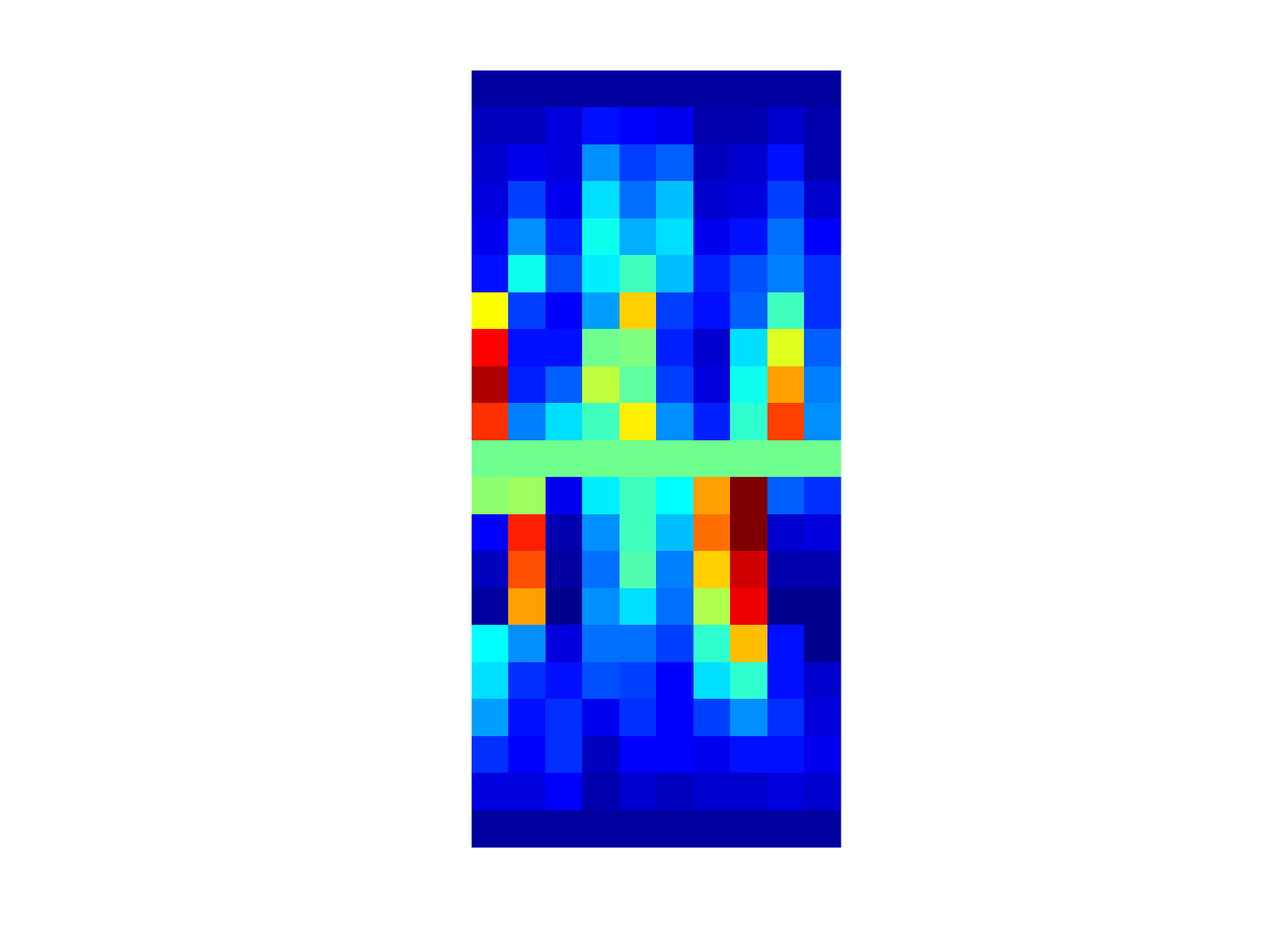}
\caption{Multiscale basis with MsFEM (x-component)}
\label{fig:basisMS} 
\end{figure}

After solving the spectral problem (\ref{spectral_1}) in $\omega_i$, we arrange the eigenvalues in  ascending order
$$0<\lambda^i_1\leq \lambda^i_2\cdots \leq \lambda^i_{J_i}.$$
Let $\phi^i_1, \phi^i_2,\cdots, \phi^i_{J_i}$ be the corresponding eigenfunctions. We select the first $l_i$ eigenfunctions whose corresponding eigenvalues less than the pre-defined tolerance to form the basis space $V^i_{\text{ms}}$,
i.e., $$V^i_{\text{ms}}=\text{span}\{\phi^i_1, \phi^i_2,\cdots, \phi^i_{l_i}\}.$$
The global velocity basis space is  $$V_{\text{ms}}=\oplus_{E_i \in \mathcal{E}_H }V^{i}_{\text{ms}}.$$
We remark that the inverse of the eigenvalue usually decays very fast. Figure \ref{fig:eig} shows the
 eigenvalue behavior corresponding to the local permeability depicted in Figure \ref{fig:localk}.
Therefore, one just needs to select a very small number of basis functions, i.e., $l_i$ is usually small, and the total dimension of $V_{\text{ms}}$   $\text{dim}(V_{\text{ms}})=\sum_i l_i$ is much less than the dimension of  $V_{\text{snap}}.$ $V_{\text{ms}}$ is the coarse space we will use in the design of the preconditioner.
We say this multiscale space the generalized multiscale space.
Figure \ref{fig:eigb1} and \ref{fig:eigb2} shows an example of the x-component of first and second spectral basis, we can see clearly that the basis includes the feature of the permeability field. 
 For comparison in numerical experiments, we also
briefly review the multiscale basis based on the mixed multiscale finite element method (MsFEM).
For each coarse neighborhood, we seek $(v_i^{ms},p_i^{ms})$ such that

\begin{equation}
\begin{split}
\kappa^{-1} v_i^{ms}+\nabla p_i^{ms}&=0,\quad\quad \text{in} ~ \omega_i,\\
\text{div}(v_i^{ms})&=c,\quad \quad \text{in}~ \omega_i.
\end{split}
\label{eq:msfem}
\end{equation}
subject to the boundary condition
\begin{equation}
\begin{split}
v_i^{ms} \cdot n_i=0,\quad \text{on} \quad\partial\omega_i,\\
v_i^{ms} \cdot n_i=1,\quad \text{on} \quad E_i
\end{split}
\end{equation}
$c$ is chosen to guarantee
the compatible condition such that the equation (\ref{eq:msfem}) is solvable.
The above problem can be solved in $K_i^1$ and $K_i^2$ respectively. This type of multiscale
basis is the same as RT0 velocity basis on the boundary of each coarse element, but
it is no longer linear (see Figure \ref{fig:basisMS}) inside the coarse element since it includes variations caused by the heterogeneity of the media.
The generalized multiscale basis can be viewed as  an enrichment of the single multiscale basis.


\section{Two-grid preconditioner}\label{sec:2grid}
In this section, we will describe the steps of  the two-grid preconditioner for solving
system (\ref{eq:reduced_fine_system}). The preconditioner generally consists
of two components,  additive local smoother and a coarse preconditioner to exchange global information.

We use the smoother introduced in \cite{pre_mixed_hdiv}, which is
the sum of some local preconditioners times a constant $\eta$.
More specifically, for each oversampled coarse block $K_i^+$ , we solve a local Neumann problem with zero flux boundary condition and the residual as
source.
The coarse preconditioner is standard, one first projects the residual into the coarse
space, and solves a coarse problem with the residual as source and
then project the solution  back to the fine-grid.

Then we can define the two-grid preconditioner as \\
{\bf Step 1:} do $m_1$ pre-smoothing steps;\\{\bf Step 2:} do one coarse correction step: transfer and solve the residual equation on the coarse level;\\
	{\bf Step 3:} do $m_2$ post-smoothing steps;
	
We can write the smoothing step $S$ in matrix form, which is:
\begin{equation}
Sr =\sum_{i=1}^{N} \eta R_i^{+T}L_i^{-1}R_i^+r
\label{eq:smooth}
\end{equation}
Where $r$ is the residual after last step.
The coarse preconditioner $P$ can also be written in matrix form as:
\begin{equation}
Pr = R_H^{T}L_H^{-1}R_Hr
\label{eq:coarsepre}
\end{equation}	
After each smoothing step or coarse correction step, we need to recompute the residual, which is the major difference between the
two-grid method and the two-level Schwarz method.
We note that the pre-smoothing  and post-smoothing  are the same and can be
completed with parallel computing easily. We use direct solver to solve all the subdomain problems and
the coarse problem.
The coarse preconditioner plays vital
role for the success of the two-grid preconditioner, using usual coarse space such as RT0 usually fails for highly heterogeneous media.
The efficiency of the preconditioner is controlled by the dimension of the coarse space, the parameters $\eta$ and
$m$, and the size of the oversampled region.

Next,  the precise algorithm to solve Equation (\ref{eq:fine_system}) can be summarized as:\\
Step 1: Perform the preprocessing step introduced in Section \ref{sec:pre} to obtain $\bar{v}_h$.\\
Step 2: Form the adaptive coarse space based on the pre-defined tolerance, and compute and store  
matrix $R_H^{T}L_H^{-1}R_H,R_i$.\\
Step 3: Use PCG combined with the two-grid preconditioner introduced above to iteratively solve Equation (\ref{eq:reduced_fine_system}).\\
Step 4: Add $\bar{v}_h$ and $\tilde{v}_h$\\
Note that step 1 and step 2 can be done independently.
\section{Numerical examples}\label{sec:numerical-results}

In this section, we present several representative numerical examples to show the performance of the two-grid preconditioner with reduced coarse space constructed by the GMsFEM as discussed above.
In all simulations reported below, the computational domain is $D = (0,1)^d, d=2,3$,
$\eta=0.2$, $K_i^+=K_i+2$, the initial guess is zero, one pre-smoothing and post-smoothing are performed. The stopping criterion for PCG iteration is that the residual is reduced by a factor of $10^7$ in $l^2$ norm if no special
declaration is specified. In all tables below, "Cond" represents the  condition number of the resulting preconditioned matrix, "$N_{iter}$" is the number of PCG iterations until pre-defined relative residual threshold is reached,
"$T_\text{setup}(s)$" is the CPU time for computing the basis and assembling the coarse matrix,
"$T_\text{solve}(s)$" is the CPU time for PCG iterations, and "Dim" is the dimension of the coarse space.
We are specifically interested in the robustness
of our method (robustness refers to the sensitivity of the convergence performance to the ratio of highest to the lowest permeability of a high contrast media) and the computational performance advantages against other two-grid preconditioners with RT0 and MsFEM space as coarse space.
All the computation is performed on a 
workstation with Intel Xeon E5-2643 CPU and with Matlab.

\begin{figure}[H]
	\centering
	\subfigure[$\kappa_{1}$]{
		\includegraphics[width=3.5in]{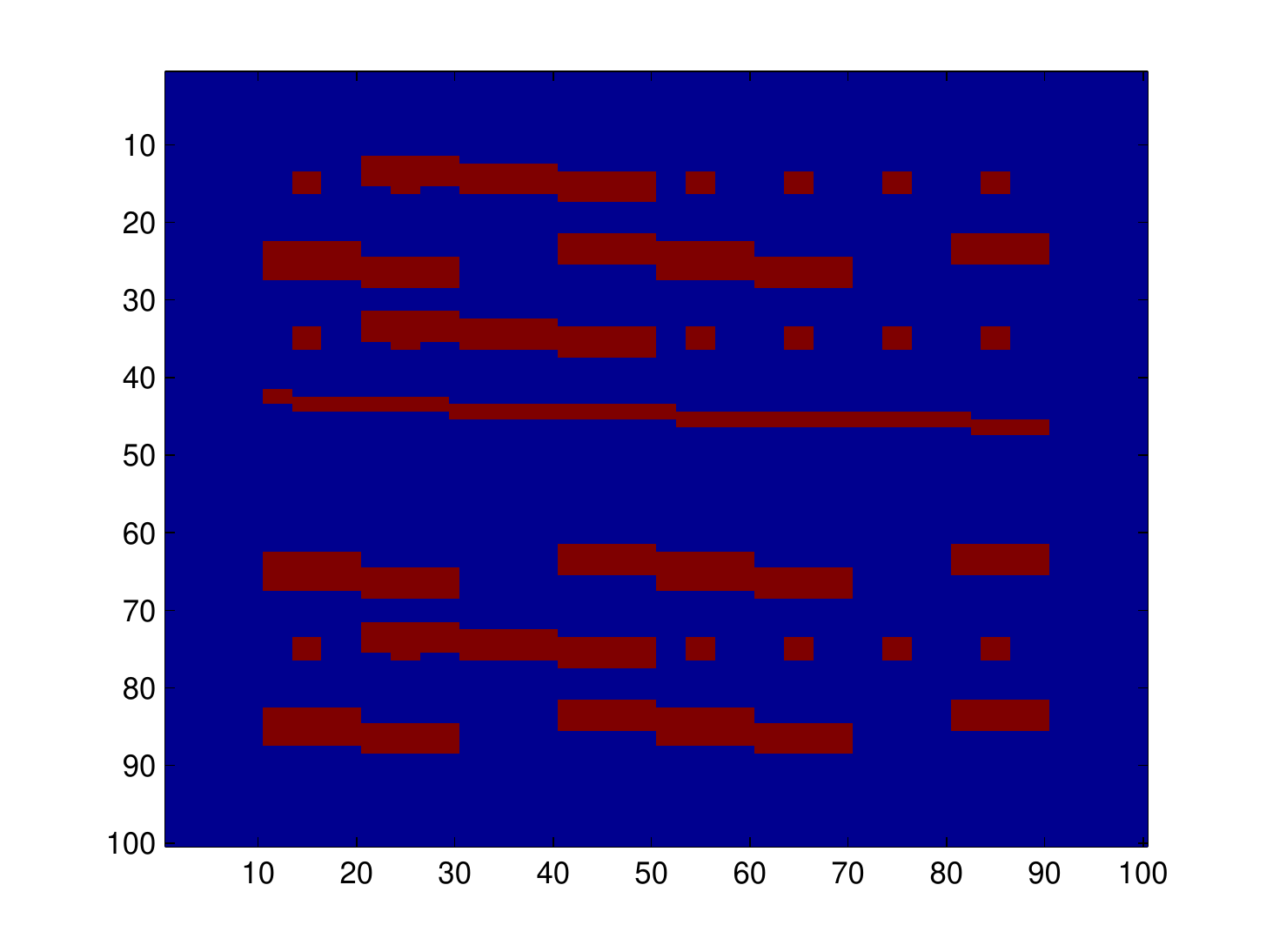}}	
	\subfigure[$\kappa_{2}$]{
		\includegraphics[width=3.5in]{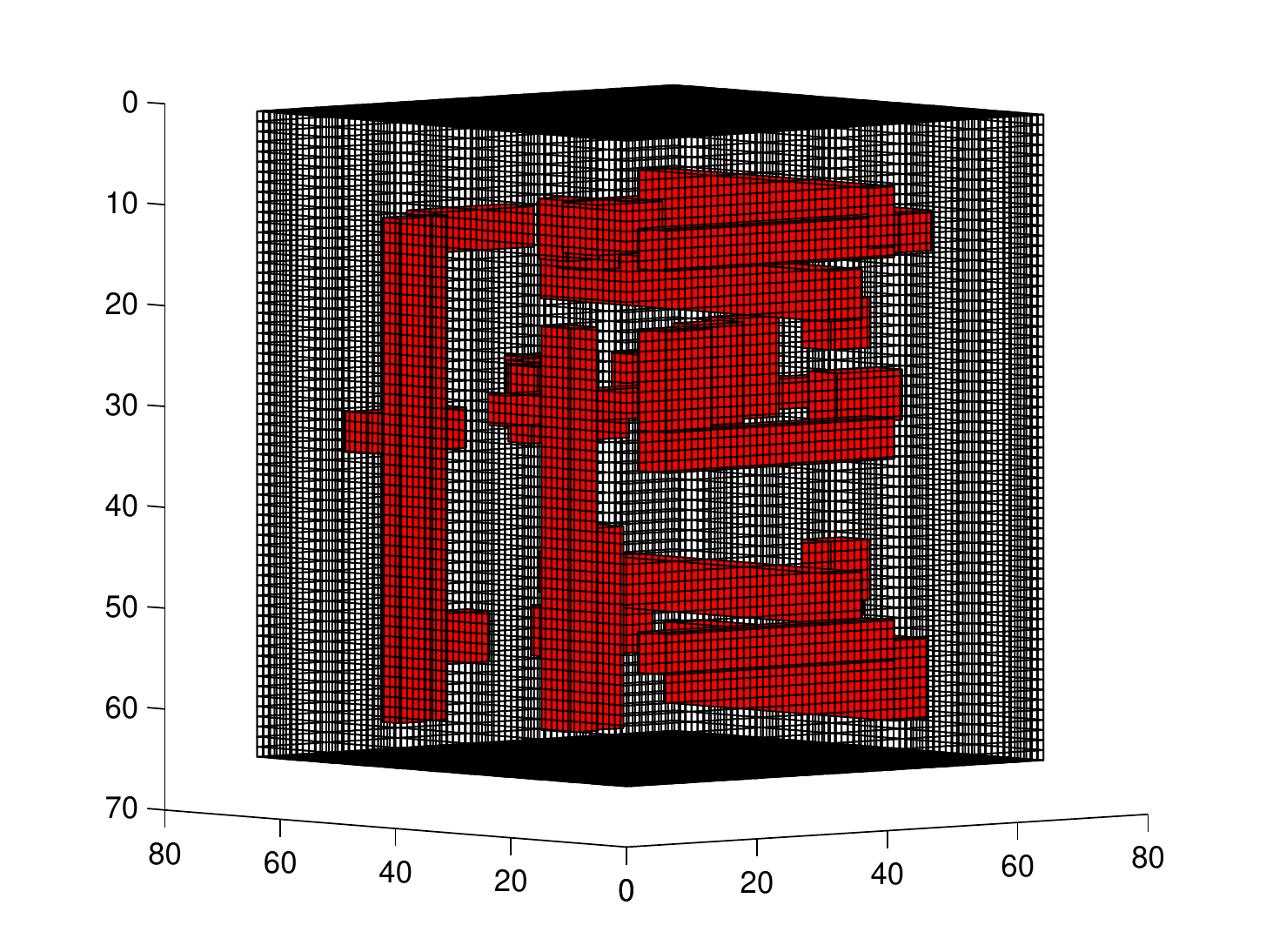}}
	\subfigure[$\kappa_{3}$ in $\log_{10}$ scale]{
		\includegraphics[width=3.5in]{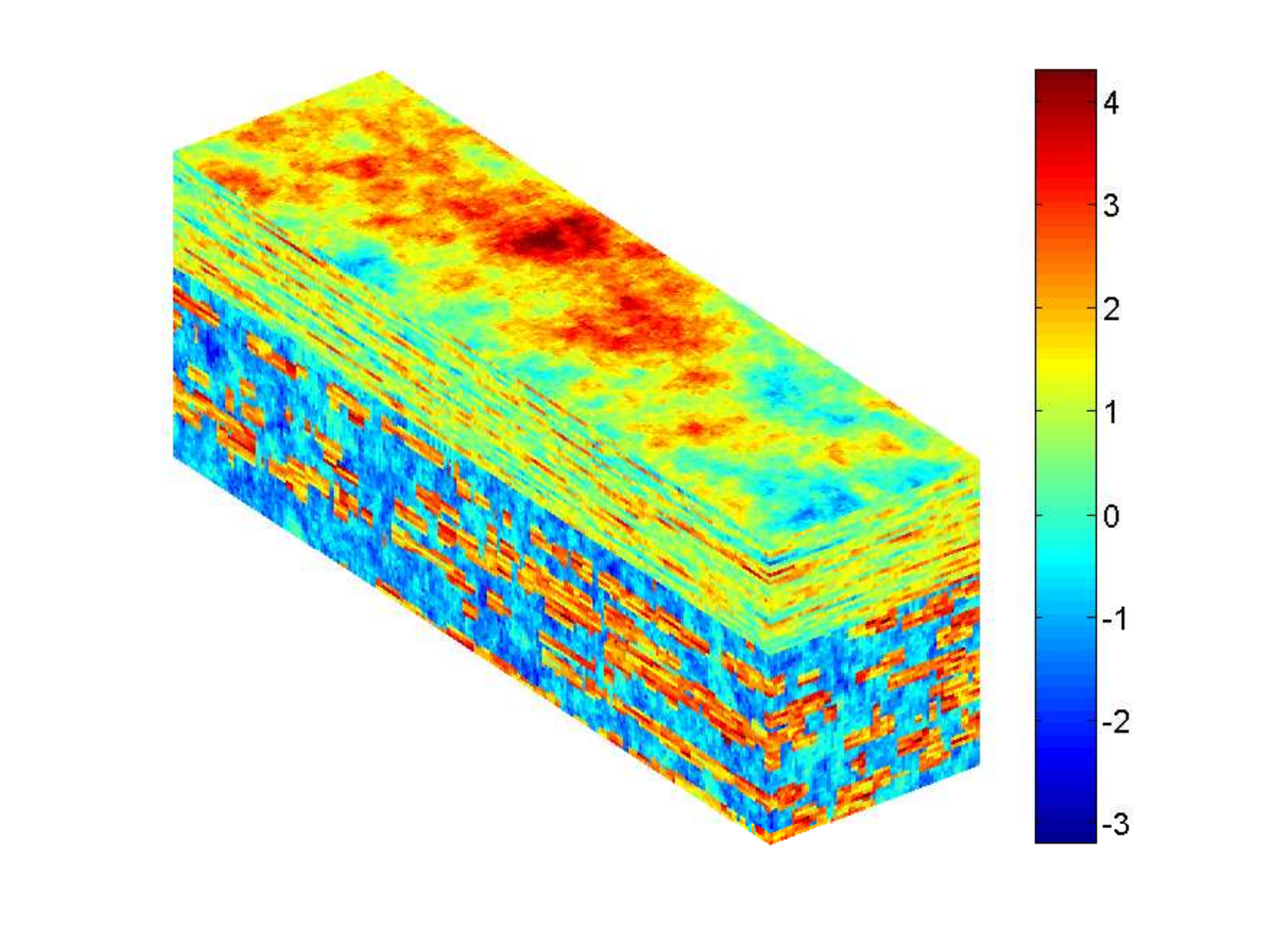}}
	\caption{Permeability fields.}
	\label{fig:models} 
\end{figure}

We consider three models  (one in 2-D and two in 3-D) with permeability $\kappa$ depicted in Figure \ref{fig:models}.  We note that
for model 1, $\kappa=1$ in the blue region and $\kappa=k_1$ in the red region, while for model 2, $\kappa=1$ in the blank region and $\kappa=k_2$ in the red region. For these two models, we will vary the orders $k_1, k_2$ to test the robustness of our method. Both models contain channels and isolated inclusions.
There are $100\times100$ fine elements and $10\times10$ coarse elements in model 1.
 For model 2, the domain $D$ is divided into $8\times 8\times 8$ square coarse elements, and in each coarse element, we generate a uniform grid with $8\times 8\times 8$ fine scale square elements.
Model 3 contains the last 80 layers of the SPE10 model \cite{Aarnes2005257}.  The precision of this model is
$220 \times 60 \times 80$, we divide the model into $22 \times 6 \times 8$ coarse blocks. The SPE 10 model is  used as a  benchmark  multiscale model in industry,  and is therefore a good test case for our methodology. We will demonstrate the computational performance of our two-grid preconditioner for two types of applications. Namely, we present Darcy flow problems and two-phase flow and transport problems. In the next subsection, we show the results of our method for solving Darcy flow problems.


\subsection{Darcy flow}

In this section, we test the computational performance and robustness of the proposed preconditioner on Darcy flow problems. First, we vary the value of $k_1$ in model 2 and $k_2$ in model 3 to show the robustness of the preconditioner, for comparison
we also present the results of preconditioner with other two coarse space,  RT0 coarse space and MsFEM coarse  space.
Tables \ref{model-1-orders_rt0}-\ref{model-1-orders} show the results for model 1. 
We can see clearly that the iteration number and condition number depend on the contrast of the media if
RT0 basis or MsFEM basis is applied, especially for the case $k_1< 0$.
However, if we use adaptive spectral coarse space, the preconditioner is  robust for both
$k_1< 0$ and $k_1\geq 0$ respectively. The dimension of the coarse space increased only 58 if 
$k_1< 0$, and there is no increase when $k_1\geq 0$. 
We observe similar phenomenon from the test results of model 2, which is reported in Tables \ref{model-2-orders_rt0}-\ref{model-2-orders}.
We therefore draw the conclusion that our preconditioner is robust against the contrast of the media.

Next, we compare the computational performance of our preconditioner with other two preconditioners that has different coarse space.  We test these three preconditioners on the two 3D models, i.e., model 2 and model 3.
The degrees of freedom of the linear system (\ref{eq:fine_system}) for these two models are 1036289 and 4188401 respectively.
 The results for model 2 and model 3 are presented  in Table \ref{model-2-methods} and Table \ref{model-3-methods} respectively.  In Table \ref{model-2-methods}, we observe that the dimension of the coarse space constructed from GMsFEM is slightly larger than the other two methods, and it takes more offline CPU time, since we need to compute snapshot space the eigenvalue problems. However, our method costs less than one half CPU time  for PCG iterations of the other two methods especially for model 3.  Moreover, the condition number of the preconditioned matrix is much less. The advantages of our method is more obvious for problems with larger size as it is shown in Table \ref{model-3-methods}, we can see that the total CPU time of our method is only about 22\% and 30\%  of the other two methods respectively. From these two examples, we can conclude that using adaptive spectral coarse space is more efficient. For larger model, we can expect more savings on computational time and memory. In the next section, we apply our method to a more realistic two-phase flow and transport problem.

\begin{table}[H]
	\centering \begin{tabular}{|c|c|c|c|c|c|c|}\hline
		$k_1$ & Dim & $N_\text{iter}$ & Cond \tabularnewline\hline
		$10^{-6}$ &282&190&2.8e+05    \tabularnewline\hline
		$10^{-4}$ &282&138&2.8e+03      \tabularnewline\hline
		$10^{-2}$ &282&43&  32.4   \tabularnewline\hline
		$10^{0}$ &282&17&3.6    \tabularnewline\hline
		$10^{2}$ &282&29&17.9     \tabularnewline\hline
		$10^{4}$ &282&31&26.3     \tabularnewline\hline
		$10^{6}$ &282&31&26.5      \tabularnewline\hline
	\end{tabular}
	\caption{Robustness test results for model 1, RT0 coarse space.}
	\label{model-1-orders_rt0}
\end{table}

\begin{table}[H]
	\centering \begin{tabular}{|c|c|c|c|c|c|c|}\hline
		$k_1$ & Dim & $N_\text{iter}$ & Cond \tabularnewline\hline
		$10^{-6}$ &282&189&1.2e+05    \tabularnewline\hline
		$10^{-4}$ &282&97&1.2e+03      \tabularnewline\hline
		$10^{-2}$ &282&33&17.2     \tabularnewline\hline
		$10^{0}$ &282&17&3.6    \tabularnewline\hline
		$10^{2}$ &282&21&7.3     \tabularnewline\hline
		$10^{4}$ &282&22&10.0     \tabularnewline\hline
		$10^{6}$ &282&22&10.0      \tabularnewline\hline
	\end{tabular}
	\caption{Robustness test results for model 1, MsFEM coarse space.}
	\label{model-1-orders_msfem}
\end{table}

\begin{table}[H]
	\centering \begin{tabular}{|c|c|c|c|c|c|c|}\hline
		$k_1$ & Dim & $N_\text{iter}$ & Cond \tabularnewline\hline
		$10^{-6}$ &340&16&3.5     \tabularnewline\hline
		$10^{-4}$ &340&16&3.5      \tabularnewline\hline
		$10^{-2}$ &340&16&3.5     \tabularnewline\hline
		$10^{0}$ &282&17&3.6      \tabularnewline\hline
		$10^{2}$ &282&18&4.1      \tabularnewline\hline
		$10^{4}$ &282&18&4.1      \tabularnewline\hline
		$10^{6}$ &282&18&4.1      \tabularnewline\hline
	\end{tabular}
	\caption{Robustness test results for model 1, GMsFEM coarse space, eigenvalue tolerance is 10.}
	\label{model-1-orders}
\end{table}


\begin{table}[H]
	\centering \begin{tabular}{|c|c|c|c|c|c|c|}\hline
		$k_2$ & Dim & $N_\text{iter}$ & Cond \tabularnewline\hline
		$10^{-6}$ &1857&28&34.8   \tabularnewline\hline
		$10^{-4}$ &1857&30&34.1   \tabularnewline\hline
		$10^{-2}$ &1857&22&12.1    \tabularnewline\hline
		$10^{0}$ &1857&15&3.1   \tabularnewline\hline
		$10^{2}$ &1857&32&17.0    \tabularnewline\hline
		$10^{4}$ &1857&40&40.1     \tabularnewline\hline
		$10^{6}$ &1857&41&40.7      \tabularnewline\hline
	\end{tabular}
	\caption{Robustness test results for model 2, RT0 coarse space.}
	\label{model-2-orders_rt0}
\end{table}
\begin{table}[H]
	\centering \begin{tabular}{|c|c|c|c|c|c|c|}\hline
		$k_2$ & Dim & $N_\text{iter}$ & Cond \tabularnewline\hline
		$10^{-6}$ &1857&18&8.9   \tabularnewline\hline
		$10^{-4}$ &1857&19&8.6    \tabularnewline\hline
		$10^{-2}$ &1857&15&4.2    \tabularnewline\hline
		$10^{0}$ &1857&15&3.1  \tabularnewline\hline
		$10^{2}$ &1857&24&9.3    \tabularnewline\hline
		$10^{4}$ &1857&32&21.2     \tabularnewline\hline
		$10^{6}$ &1857&32&21.5      \tabularnewline\hline
	\end{tabular}
	\caption{Robustness test results for model 2, MsFEM coarse space.}
	\label{model-2-orders_msfem}
\end{table}

\begin{table}[H]
\centering \begin{tabular}{|c|c|c|c|c|c|c|}\hline
$k_2$ & Dim & $N_\text{iter}$ & Cond \tabularnewline\hline
$10^{-6}$ &1904&12&3.1     \tabularnewline\hline
$10^{-4}$ &1904&12&3.1      \tabularnewline\hline
$10^{-2}$ &1898&12&3.0      \tabularnewline\hline
$10^{0}$ &1857&13&3.0      \tabularnewline\hline
$10^{2}$ &1912&14&3.2      \tabularnewline\hline
$10^{4}$ &1922&14&3.3      \tabularnewline\hline
$10^{6}$ &1941&14&3.4      \tabularnewline\hline				
\end{tabular}
\caption{Robustness test results for model 2, GMsFEM coarse space, eigenvalue tolerance is 10.}
\label{model-2-orders}
\end{table}

\begin{table}[H]
\centering \begin{tabular}{|c|c|c|c|c|c|c|}\hline
Coarse space &Dim&$T_\text{setup}(s)$   & $T_\text{solve}(s)$    & $N_\text{iter}$&Cond    \tabularnewline\hline
RT0&1857&2.1&58.9& 38& 40.1    \tabularnewline\hline
MsFEM&1857&5.4&45.8&30 &21.2       \tabularnewline\hline			
GMsFEM&1922& 20.9&21.7&14&3.3    \tabularnewline\hline
\end{tabular}
\caption{Computational performance tests for model 2 with different coarse space, $k_2=10^4$, eigenvalue tolerance is 10.}
\label{model-2-methods}
\end{table}

\begin{table}[H]
\centering \begin{tabular}{|c|c|c|c|c|c|c|}\hline
Coarse space &Dim&$T_\text{setup}(s)$   & $T_\text{solve}(s)$    & $N_\text{iter}$&Cond    \tabularnewline\hline
RT0&3869&11.1&1243.4&196 & 847.7    \tabularnewline\hline
MsFEM&3869&49.2&856.8&135 & 397.2       \tabularnewline\hline			
GMsFEM&7320&160.2&122.4 &18 & 9.2  \tabularnewline\hline
\end{tabular}
\caption{Computational performance test results for model 3 with different coarse space, eigenvalue tolerance is 10.}
\label{model-3-methods}
\end{table}

\subsection{A two phase flow and transport problem}
In this section, we test our method for the case where two immiscible fluid phases, i.e., water and oil, are flowing in a heterogeneous porous media. In particular, we consider two-phase flow in a reservoir domain (denoted by $\Omega$).
First, we summarize the underlying partial differential equations \cite{efendiev2016online,SpeJ}.  The basic equation describing the filtration of a fluid through a porous media is the continuity equation, which states that mass is conserved (assuming that the rock and fluids are incompressible):
\begin{equation} \label{mass}
\phi\frac{\partial s_l}{\partial t} + \nabla \cdot {\bf v}_l = q_l, \quad l=o, w.
\end{equation}
where ${  \bf v}_l$ is the phase velocity for phase $l$, $s_l$ is
saturation, $o$ and $w$ refers to the oil and water phases, $\phi$ is the porosity of the medium, and $q_l$ is the source term, which models sources and sinks, i.e., outflow and inflow per volume at designated well locations.

If ignore both gravity and capillary pressure effects, for each phase, phase velocity is related to pressure $p$ by the Darcy's law:
\begin{equation} \label{darcy}
{ \bf v}_l=-{ K}\frac{k_{rl}(s_l)}{\mu_l}  \nabla {p}
\end{equation}
$ {K}$ is the absolute permeability
tensor, $k_{rl}$ is the relative permeability to phase $l$ ($l=o, w$), $\mu_j$ is the viscosity,  and ${p}$ is pressure. We denote $\lambda_l=\frac{k_{rl}(s_l)}{\mu_l} $ as the phase mobility.
The relative permeability for a phase is usually a nonlinear function of the saturation of that phase. Throughout the paper, we use a single set of relative permeability.

Combining Darcy's law, mass conservation, and the property $s_w+s_o=1$, we derive the following coupled system of pressure and saturation equations (we use $s$ instead of $s_w$ for simplicity):
\begin{eqnarray}
\nabla\cdot {\bf v} &=& q_w + q_o \quad \textrm{in} \quad \Omega \label{pressEq1}\\
\phi \frac{\partial s}{\partial t}+ \nabla \cdot({ f_w(s) {\bf v}}) &=& \frac{q_w}{\rho_w} \quad \textrm{in} \quad \Omega \label{satEq}\\
{\bf v}\cdot n & = & 0 \quad \textrm{on} \quad \partial {\Omega} \quad \textrm{(no flow at boundary)}\\
s(t=0) & = & s_{0}, \quad \textrm{in} \quad \Omega \quad \textrm{(initial known saturation)} \label{initial}
\end{eqnarray}

where ${\bf v} = { \bf v}_w + { \bf v}_o$ is the total velocity, which can be expressed as
\begin{equation} \label{pressure equ}
{\bf v} =-\lambda(s) { K}  \nabla {p}.
\end{equation}
Equation (\ref{pressEq1}) is the overall continuity equation, which is referred as the "pressure equation".

$f_w(s)$ is the flux function,
\begin{equation} \label{fraction}
f_w(s) =\frac{\lambda_{w}(s)}{\lambda(s)}=\frac{k_{rw}(s)}{k_{rw}(s)+\frac{\mu_w}{\mu_o}k_{ro}(s)}
\end{equation}

 where $\lambda$ is the total mobility defined as

\begin{equation} \label{mobility}
\lambda(s) =\lambda_w(s)+\lambda_o(s)=\frac{k_{rw}(s)}{\mu_w}+\frac{k_{ro}(s)}{\mu_o}
\end{equation}
Equation (\ref{satEq}) is referred as the "saturation equation", obtained by writing the continuity equation for phase $w$ in terms  of the total velocity.

Equations (\ref{pressEq1})-(\ref{initial}) is a nonlinear coupled system.  The coupling  is through the saturation-dependent mobilities $\lambda_l$ in the pressure equation, and the pressure-dependent velocity in the saturation equation.The solution strategy we  use for this system is a  sequential splitting method, called the IMPES, which is widely used in reservoir simulation. IMPES treats the flow and transport separately and differently. That is at each time step one solves for the pressure and velocity first and then uses the velocity to solve for the saturation. Note that the total mobility $
\lambda$ in Equation (\ref{pressure equ}) is calculated at the previous time level.

After obtaining the saturation distribution at the end of a time step, the total mobility $\lambda$ is updated, and a system of new pressure equations is assembled.The calculation of the global pressure solution on the underlying fine grid at each time step is the most time consuming part of the strategy.   The involving challenge can be observed from Equation (\ref{pressure equ}). The total mobility $\lambda$  is a function of space and time. The absolute permeability $ {K}$  is usually the dominant part in dictating the flow field in natural porous formations. The variability and correlation structure of permeability are usually  expressed as a complex multiscale function. Moreover, $ {K}$  often displays significant variation within small distance. Therefore, capturing the variability and resolving the spatial correlation structures usually result in heavy computational burden. 

We note that the two-grid preconditioner  for the mixed formulation of the pressure equation is particularly suitable for the flow and transport problem, since it can compute pressure and velocity simultaneously. The saturation equation is solved by finite volume method, together with backward Euler scheme for time discretization. Specifically, consider a fine grid cell
$\Omega_{i}$ with edges $\gamma_{ij}$ and associated normal vectors
$n_{ij}$ pointing out of $\Omega_{i}$, the saturation Equation (\ref{satEq}) is discretized as

\begin{equation}\label{eq:DisSatEq_Implicit}
 s_i^{n+1}= s_i^n +\frac{\Delta t}{|\Omega_i|}\left(q^+-\sum_j F_{ij}(s^{n+1})u_{ij}+f_w(s_i^{n+1})q^-\right).
\end{equation}
 where $s_{i}^{n}$ is the cell-average of the water saturation at time
$t=t_{n}$, $q^+ = \max(q_i,0)$ and $q^-=\min(q_i,0)$,  ${u}_{ij}$ is the total velocity (for oil and water) over the edge $\gamma_{ij}$ between the two adjacent cells and $F_{ij}$ is a numerical approximation of the flux over edge $\gamma_{ij}$ defined as,
\begin{equation}
F_{ij}\approx\int_{\gamma_{ij}}\left(f_{ij}(s)v_{ij}\right).n_{ij} \; dv. \label{eq:integralFlux}
\end{equation}
There are different schemes to evaluate the integrand in Equation (\ref{eq:integralFlux}). A common approach is to use a first order approximation, known as upstream weighting that is defined as,
 \begin{equation}
 f_{ij}(s)=\begin{cases}
   f_w(s_i)      &\text{if $\; v_{ij}\cdot n_{ij}\geq0$;}\\
   f_w(s_j)      &\text{if $\; v_{ij}\cdot n_{ij}< 0$.}
 \end{cases}
 \end{equation}
The saturation equation  Equation (\ref{eq:DisSatEq_Implicit}) yields a  nonlinear system and can be solved for $s^{n+1}$ ( $s^{n+1}$ is the vector of cell-saturations at time instance $t_{n+1}$) by iterative methods, such as Newton-Raphson,  efficiently. 

\begin{figure}[H]
	\centering
	\subfigure[Saturation at $t=50$]{
		\includegraphics[width=3.5in]{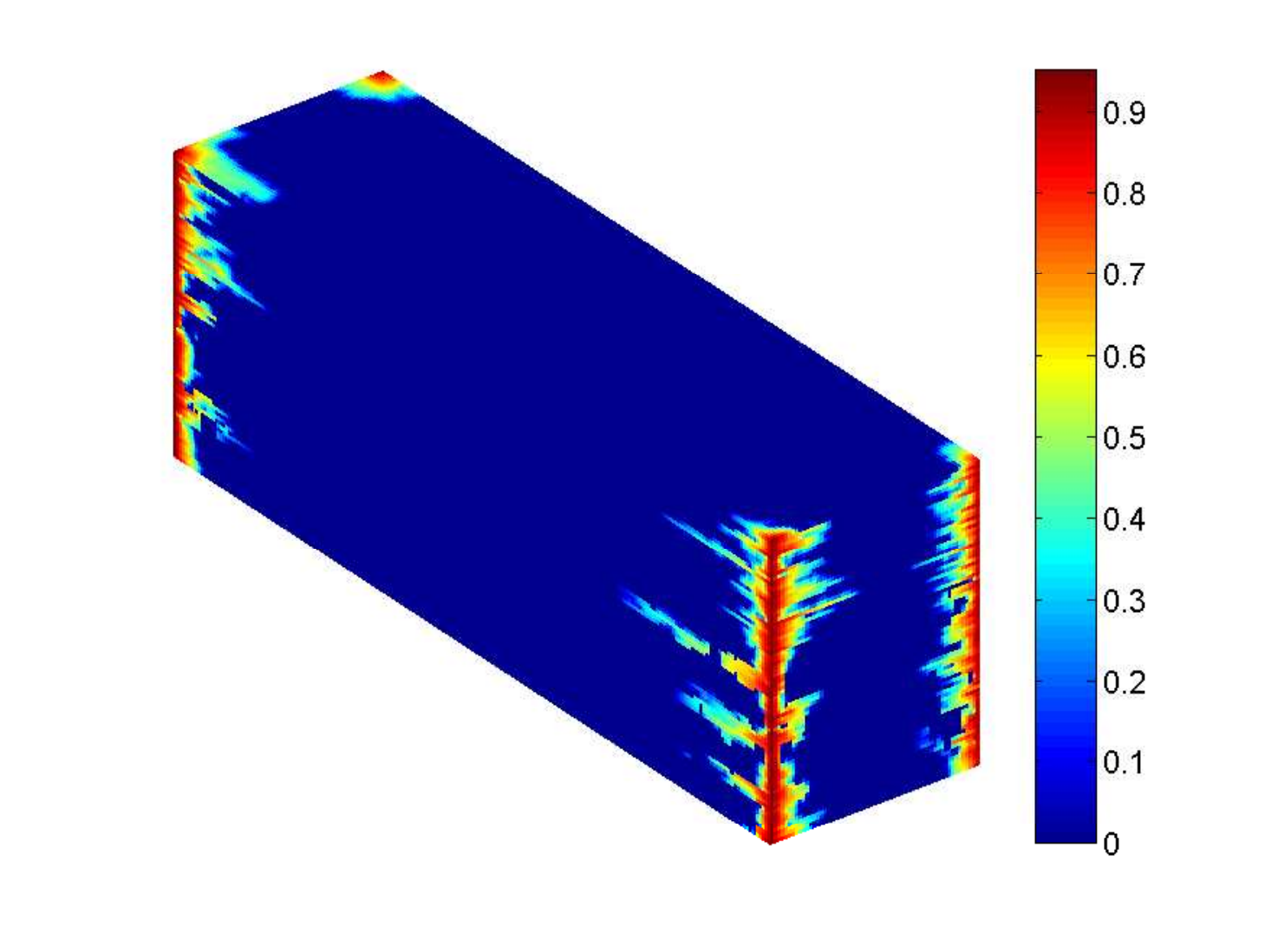}}		
	\subfigure[Saturation at $t=1000$]{
		\includegraphics[width=3.5in]{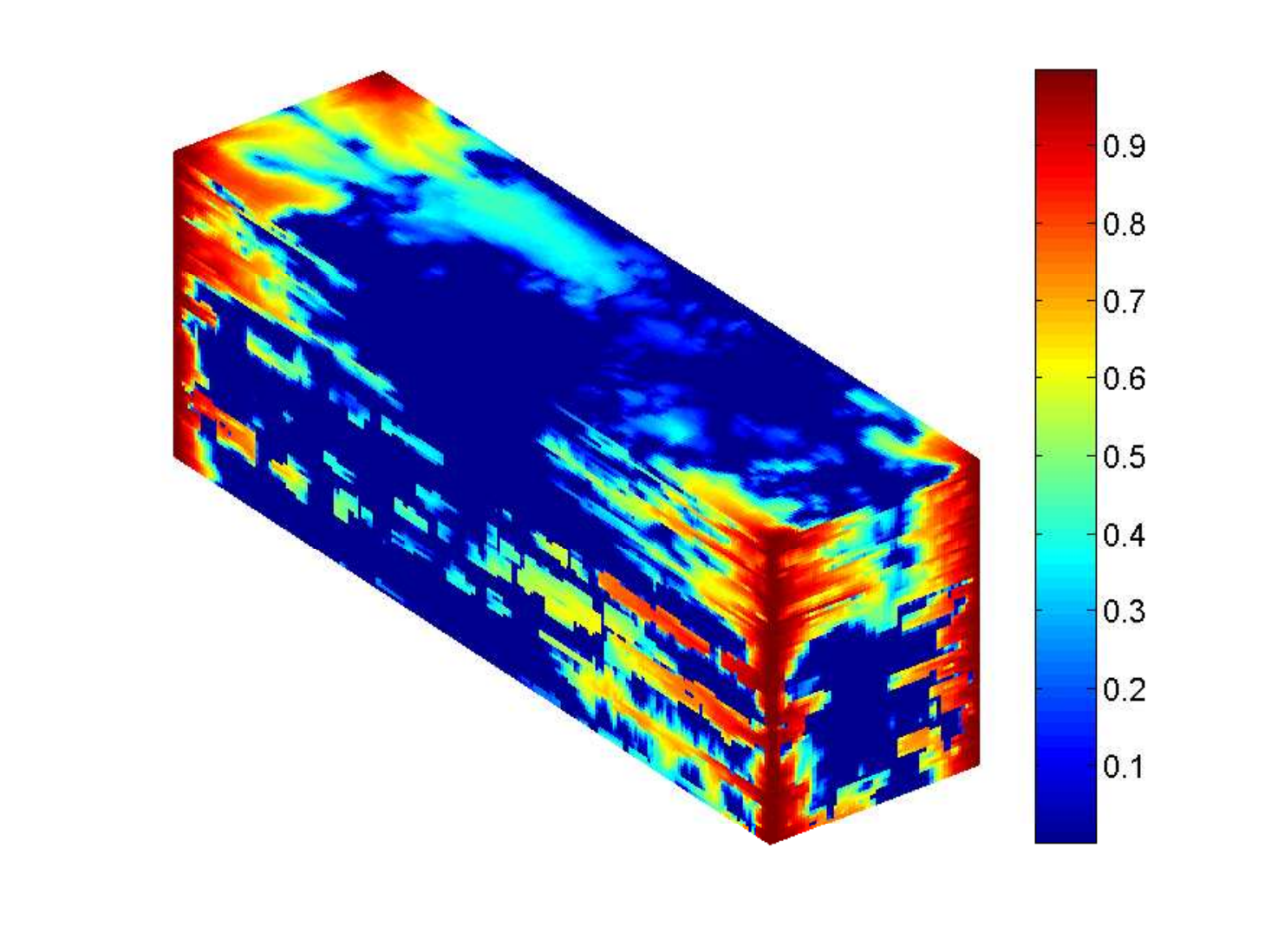}}	
	\subfigure[Saturation at $t=2000$]{
		\includegraphics[width=3.5in]{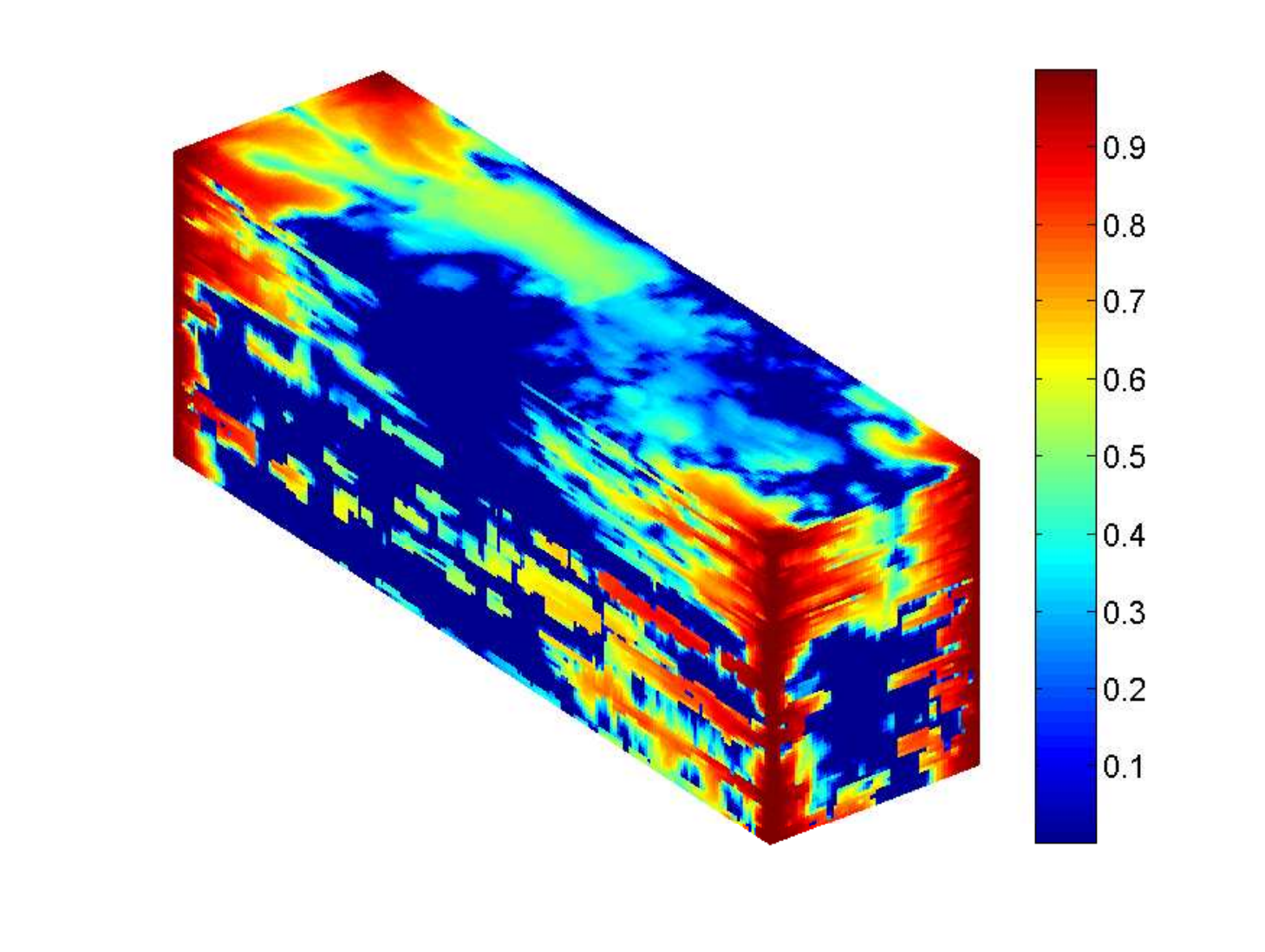}}		
	\caption{Saturation profiles for model 3}
	\label{fig:Sat_model2}
\end{figure}

During the whole simulation process, we need to solve a number of times of the pressure equation with varying coefficients.
However we don't need to update the adaptive coarse space. Instead, we use the multiscale space calculated at the initial time step. We test our method on model 3, that is,  the last 80 layers of the SPE10 model. The initial saturation field is zero, i.e., the reservoir is filled with oil at initial time.  The viscosity for water and oil are 1 and 5 respectively. The time steps simulated is 2000 in total, and we  solve the pressure equation 40 times. The eigenvalue tolerance is 15, and corresponding dimension of the coarse space is 9941. We are interested in solving the pressure equation with different accuracy requirements via our preconditioner. Table \ref{model-3-comp time} presents the computation time of direct solver and our iterative solver for the two phase flow and transport problem. It takes 64.66 hours for the direct solver. While for the iterative solver, by setting the relative residual tolerance as $10^{-2}, 10^{-4}, 10^{-7}$, the CPU time is 13.56 hours, 1.58 hours, and 2.05 hours respectively, which are all much less than the time of the direct solver. The reason why the case $10^{-2}$, takes more time
than the case $10^{-4}$ and $10^{-7}$ is that the accuracy of solving the pressure equation has huge
influence on the computational time of solving the transport equation. 
If the pressure equation is not solved very accurately, then solving \ref{eq:DisSatEq_Implicit} will be 
expensive. Clearly, if we set the residual tolerance as $10^{-2}$, the accuracy of the velocity 
is not enough.
The computation time for the case $10^{-7}$ is 2.05 hours, which is greater than the case of $10^{-4}$.  The explanation for 
this is that, the accuracy for solving the pressure equation is sufficient for both relative residual tolerances, therefore, now the dominant factor is not the computational time for the saturation equation, but the computational time for solving all 40 times of the pressure equation.
We also want to mention that although the coarse space is fixed, the iteration
number for solving the pressure equation in later time instant will only increase 1 or 2
compared with using the coarse space computed from the exact permeability.
 The saturation profiles  at time 50, 1000 and 2000 are depicted in Figure \ref{fig:Sat_model2}. We can see with the advancement of time, water is injected into the wells in the four corners, driving oil flowing toward the production well in the middle. The moving profile of water is influenced by the multiscale features of the background permeability field.

\begin{table}[H]
	\centering 	\begin{tabular}{|c|c|c|c|c|c|c|}\hline
		Method  & relative residual & $T(h)$     \tabularnewline	\hline
		direct solver&0 & 64.66     \tabularnewline\hline
		iterative solver&$10^{-2}$ & 13.56 \tabularnewline\hline
		iterative solver&$10^{-4}$	& 1.58    \tabularnewline\hline 
		iterative solver&$10^{-7}$	&2.05    \tabularnewline\hline
	\end{tabular}
\caption{Comparison of computational time  for model 3 (two-phase flow simulation), eigenvalue tolerance is 15, dimension of the coarse system is 9941.}
\label{model-3-comp time}
\end{table}

\section{Conclusions}
In this paper, we propose a two-grid preconditioner  for mixed formulations of elliptic problems in highly heterogeneous  porous media. The main novelty  is that we use the multiscale space constructed from the GMsFEM for the coarse preconditioner. The multiscale space consists of basis functions that can capture the multiscale feature of the underlying permeability field. Preprocessing steps are used to transform the indefinite saddle problem to a positive definite one. We present numerical results to show that our preconditioner is robust in terms of contrast orders of the permeability. By comparing to other preconditioners that incorporate RT0 and the standard MsFEM space as coarse preconditioner, we demonstrate that our preconditioner is more efficient and robust. Moreover, we apply our method to a more realistic incompressible two-phase flow and transport problem, and results show that the method is highly efficient.

\section*{Acknowledgements}

EC's work is partially supported by Hong Kong RGC General Research Fund (Project 14304217)
and CUHK Direct Grant for Research 2017-18.

\bibliographystyle{plain}
\bibliography{references}

\end{document}